\documentstyle[12pt,amscd,amssymb]{amsart}
\theoremstyle{plain}
 \newtheorem{thm}{Theorem}[section]
 \newtheorem{lem}[thm]{Lemma}
 \newtheorem{prop}[thm]{Proposition}
 \newtheorem{cor}[thm]{Corollary}
 \newtheorem{guess}[thm]{Conjecture}
\theoremstyle{definition}
 \newtheorem{defn}{Definition}[section]
\theoremstyle{remark}
 \newtheorem{rem}{Remark}[section]

\newcommand{\Ext}{\operatorname{Ext}}
\newcommand{\Hom}{\operatorname{Hom}}

\newcommand{\im}{\operatorname{im}}

\newcommand{\rk}{\operatorname{rk}}

\font\b=cmr10 scaled \magstep5
\def\bigzerou{\smash{\lower1.7ex\hbox{\b 0}}}

\numberwithin{equation}{section}

\begin{document}

\title
{Euler characteristics of $SU(2)$ instanton moduli spaces
on rational elliptic surfaces}
\author{K\={o}ta Yoshioka\\
 Department of Mathematics\\
Faculty of Science, Kobe University}
 \address{Department of Mathematics, Faculty of Science,
Kobe University,
Kobe, 657, Japan }
\email{yoshioka@@math.kobe-u.ac.jp}
 \maketitle
\pagestyle{plain}

\section{Introduction}
In [MNVW], Minahan, Nemeschansky,
Vafa and Warner computed ``Euler characteristics''
of $SU(2)$-instanton moduli spaces on rational elliptic surfaces.
In order to state their formula,
let us introduce some notations.  
Let $\pi:X \to {\Bbb P}^1$ be a rational elliptic surface,
$s$ a section of $\pi$ and $f$ a fiber of $\pi$.
For $x,y \in H^2(X,{\Bbb Z})$,
we define a pairing $(x,y)_*:=-(x,y)$,
where $(x,y)$ is the intersection pairing.
Then the orthogonal complement 
$$
\langle s,f \rangle^{\perp}:=
\{x \in H^2(X,{\Bbb Z})| (x,s)_*=(x,f)_*=0 \}
$$
is isomorphic to the $E_8$-lattice.
Under the action of translations by $2 x, x \in E_8$ 
and the Weyl group $W(E_8)$,
there are three orbits which are represented by
$v_0, v_{even}, v_{odd}$ respectively,
where $v_0=0$, $(v_{even}^2)_*=4$, and 
$(v_{odd}^2)_*=2$.
We shall choose a Hodge metric $g$ such that the size of $f$
is very small.
Let $V(v_{\lambda},\Delta)$ be the instanton moduli spaces
of topological invariants $(v_{\lambda},\Delta)$.
Let $Z_{v}(\tau), v \in H^2(X,{\Bbb Z})$ be the partition function
considered by Vafa and Witten in [V-W].
Up to holomorphic anomaly, they are generating functions
of $``\chi(V(v,\Delta))"$: 
\begin{equation}
Z_{v}(\tau)=q^{-2\chi(X)/24}
\sum_{\Delta}``\chi(V(v,\Delta))"q^{\Delta},
\end{equation}
where $q=\exp(2 \pi \sqrt{-1} \tau)$.

Recently,
Minahan, Nemeschansky, Vafa and Warner
computed related partition function $Z_2$
( numbers of BPS states) [MNVW, (6.14), (6.15)].
Since the number of BPS states is minus the Euler characteristics,
their formulas (6.15) are stated as follows:
\begin{align*}
Z_0(\tau)&=\frac{-1}{24 \eta(\tau)^{24}}
\left(\widehat{E}_2(\tau)P_{0}(\tau)+
\left(\theta_3(\tau)^4\theta_4(\tau)^4-
\frac{\theta_2(\tau)^8}{8}\right)
(\theta_3(\tau)^4+\theta_4(\tau)^4)\right),\\
Z_{even}(\tau)&=\frac{-1}{24 \eta(\tau)^{24}}
\left(\widehat{E}_2(\tau)\frac{P_{even}(\tau)}{135}-
\theta_2(\tau)^8 \frac{(\theta_3(\tau)^4+\theta_4(\tau)^4)}{8}\right),\\
Z_{odd}(\tau)&= \frac{-1}{24\eta(\tau)^{24}}
\left(\widehat{E}_2(\tau) \frac{P_{odd}(\tau)}{120}
-\frac{1}{8}\theta_2(\tau)^4 E_4(\tau)
\right),
\end{align*}
where $Z_{\lambda}(\tau)=(Z_{v_{\lambda}}(\tau)+Z_{f+v_{\lambda}}(\tau))/2$.
For the definition of $P_{\lambda}(\tau)$, see \eqref{eq:p}.

On the other hand,
it is expected that there is a good compactification 
$\overline{V}(v_{\lambda},\Delta)$ of
$V(v_{\lambda},\Delta)$ such that
$``\chi(V(v_{\lambda},\Delta))"$ can be thought of some Euler
characteristics of $\overline{V}(v_{\lambda},\Delta)$.
In particular, we may expect that
if the Gieseker's compactification $M(v_{\lambda},\Delta)$
is smooth
(in our case, $v_{\lambda} \ne 0$ or $\Delta \not\in 2{\Bbb Z}$), 
then $``\chi(V(v_{\lambda},\Delta))"$ coincides with 
the ordinary Euler characteristics $\chi(M(v_{\lambda},\Delta))$.

In this paper, 
we shall compute $\chi(M(v_{\lambda},\Delta))$ and compare two
results (Theorem \ref{thm:1}).
In particular, we show that
they coincide up to holomorphic anomaly coming from
$E_2(\tau)$ and unknown terms coming from
singular spaces!

In section 1, we collect some equalities on (quasi) modular forms.
In section 2, we shall consider moduli of stable sheaves on rational 
elliptic surfaces.
In particular, we shall define our invariant 
$\widetilde{Z}_{\lambda}(\tau)$ and get a description of Poincar\'{e}
polynomials of moduli spaces.
Since $X$ is 8 points blow-ups of $\Sigma_1$,
there is a morphism $\psi:X \to {\Bbb P}^1$ such that 
general fibers are ${\Bbb P}^1$.
Let $g \in  H^2(X,{\Bbb Z})$ be the cohomology class of a fiber.
Then we get $(f^2)=(g^2)=0$.
Let $\Theta^{f,g}_{H^2(X,{\Bbb Z})}$ be the Jacobi form of 
G\"{o}ttsche and Zagier [G-Z1] associated to $f,g$.
Then Poincar\'{e}
polynomials of moduli spaces are related to this Jacobi form.  
Since the orthogonal complement of $\langle f,g \rangle$
is the $D_8$-lattice,
we can describe our invariants by using $D_8$-theta functions. 
$E_2(\tau)$ also appears as a contribution of $\langle f,g \rangle$.
These calculations are treated in section 3. 

After wrote up the first version of this paper,
the author heard that G\"{o}ttsche and Zagier [G-Z2] related 
Hodge numbers  
of moduli spaces to Jacobi forms
in a more general situation.
In particular, they showed that Euler characteristics and
signatures have generating functions in terms of 
(quasi) modular forms.
If moduli space is smooth and has a universal family, 
then Beauville [B] showed that cohomologies classes are generated by
algebraic cycles.
Hence Poincar\'{e} polynomials in section 3.3 also give
presentations of Hodge numbers of moduli spaces.   
\vspace{1pc}

{\it Notation.}

Let $X$ be a smooth projective surface.
For a torsion free sheaf $E$ on $X$,
we set 
$$
\Delta(E):=c_2(E)-\frac{\rk(E)-1}{2 \rk(E)}(c_1(E)^2)
 \in {\Bbb Q}.
$$
Let $H$ be an ample divisor on $X$.
For a pair $(c_1,\Delta) \in H^2(X,{\Bbb Z}) \times {\Bbb Q}$,
let $M_H(c_1,\Delta)$ be the moduli space of semi-stable
sheaves $E$ of rank 2 on $X$ such that
$c_1(E)=c_1$ and $\Delta(E)=\Delta$.

For a rational surface $X$,
we set
$$
Z_t(X,u):=
\frac{1}{(1-u)^{b_0(X)}(1-tu)^{b_2(X)}(1-t^2u)^{b_4(X)}}.
$$
If $X$ is defined over a finite field ${\Bbb F}_t$
of $t$ elements,
then $Z_t(X,u)$ is Weil's zeta function of $X$.

For a smooth manifold $M$,
$$
P(M,t):=\sum_i b_i(M) t^{i/2}
$$
is the Poincar\'{e} polynomial of $M$.

We set
\begin{alignat*}{2}
\theta_2(\tau) &:= \sum_{n \in {\Bbb Z}+1/2}q^{n^2/2} &\quad
E_2(\tau)&:=1-24\sum_{n >0} \sigma_1(n) q^n\\
\theta_3(\tau) &:= \sum_{n \in {\Bbb Z}}q^{n^2/2} &\quad
E_4(\tau)&:=1+240\sum_{n >0} \sigma_3(n) q^n\\
\theta_4(\tau) &:= \sum_{n \in {\Bbb Z}}(-1)^n q^{n^2/2}& \quad
e_1(\tau) &:=-\frac{1}{6}(1+24( \sum_n \sigma_1^{odd}(n)q^n))\\
\Theta(\tau)&:=\theta_3(2\tau) &\quad
F(\tau) &:=\sum_{n>0, 2 \nmid n}\sigma_1(n)q^n,
\end{alignat*}
where $q=\exp(2 \pi \sqrt{-1} \tau)$,
$\sigma_k(n)=\sum_{d|n}d^k$, $k=1,3$ and
 $\sigma_1^{odd}(n)=\sum_{d|n, 2 \nmid d}d$.
Since $E_2(\tau)$ is not modular,
we also consider $\widehat{E}_2(\tau):=E_2(\tau)-3/(\pi \im \tau)$.

For a positive definite lattice $L$ with a pairing $(\quad,\quad)$,
we define the lattice theta function:  
\begin{equation}
\Theta_{L}(x,u):=\sum_{n \in L} u^{(n^2)/2}e^{(n,x)},
\end{equation}
where $u:=\exp(2 \pi \sqrt{-1} \tau)$ and $e:=\exp(2 \pi \sqrt{-1})$.
For $\lambda \in L_{\Bbb R}$,
we set
\begin{equation}
(\Theta_{L}|\lambda)(x,u)
:=\sum_{n \in L+\lambda/2} u^{(n^2)/2}e^{(n,x)}.
\end{equation}

\begin{rem}
In the usual notation [G-Z1],
 $\Theta_{L}|\lambda=\Theta_{L}|(\lambda,0)$.
\end{rem}

\section{Preliminaries}
We shall collect some equalities which will be used later.
\begin{lem}\label{lem:wt2}
\begin{align*}
e_1(\tau) & =-\frac{1}{6}(-E_2(\tau)+2 E_2(2 \tau)),\\
F(\tau) &=- \frac{1}{24}(E_2(\tau)-3 E_2(2 \tau)+2 E_2(4 \tau)),\\
\sum_{d | n, 2 \nmid (n/d)}
d q^n &=\frac{E_2(2\tau)-E_2(\tau)}{24},\\
\sum_{d | n, 2 \nmid (n/d)}
(-1)^d d q^n &=\frac{E_2(\tau)-5 E_2(2 \tau)+4 E_2(4 \tau)}{24}. 
\end{align*}
\end{lem}

\begin{lem}\label{lem:e1}
\begin{equation}
\begin{split}
-6 e_1(\tau) & =\frac{\theta_3(\tau)^4+\theta_4( \tau)^4}{2}
=\Theta(\tau)^4+16 F(\tau),\\
\theta_4(2\tau)^4 & =\Theta(\tau)^4-16 F(\tau),\\
\theta_2(2 \tau)^4 & =16 F(\tau),\\
\theta_2(\tau)^8 &= 16 \Theta(\tau)^4 16 F(\tau).
\end{split}
\end{equation}
\end{lem}

\begin{pf}
We note that $\theta_i(2 \tau)^4\; (i=2,3,4), \theta_2(\tau)^8,
\theta_3(\tau)^4+\theta_4(\tau)^4,
e_1(\tau)$ and $F(\tau)$ are modular form for $\Gamma_0(4)$.
Hence these are polynomials of $\Theta(\tau)$
and $F(\tau)$.
Comparing Fourier coefficients of small orders,
we get our lemma.
\end{pf}

\begin{cor}\label{cor:t1}
$$
\theta_2(\tau)^4=\theta_3(\tau)^4-\theta_4(\tau)^4.
$$
\end{cor}
It is known that
\begin{align*}
\Theta_{E_8}(0,\tau) &=
 \frac{\theta_2(\tau)^8+\theta_3(\tau)^8+\theta_4(\tau)^8}{2}\\
&=E_4(\tau).
\end{align*}
We set
\begin{equation}\label{eq:p}
\begin{split}
P_0(\tau)&:=E_4(2 \tau),\\
P_{even}(\tau)&:=\frac{E_4(\tau/2)+E_4(\tau/2+1/2)}{2}-E_4(2\tau),\\
P_{odd}(\tau)&:=\frac{E_4(\tau/2)-E_4(\tau/2+1/2)}{2}.
\end{split}
\end{equation}
Then we obtain the following lemma.

\begin{lem}\label{lem:wyle}
\begin{align}
P_0(\tau) &=\frac{\theta_2(2 \tau)^8+ \theta_3(2 \tau)^8+
\theta_4(2 \tau)^8}{2},\\
\frac{P_{even}(\tau)}{135}&=
\frac{\theta_2(2 \tau)^8+ \theta_3(2 \tau)^8-
\theta_4(2 \tau)^8}{2},\\
\frac{P_{odd}(\tau)}{120} &=
\frac{\theta_2(2 \tau)^6 \theta_3(2 \tau)^2+ \theta_2(2 \tau)^2
\theta_3(2 \tau)^6}{2}.
\end{align}
\end{lem}

\begin{pf}
The first equality is trivial.
We note that $P_{even}(\tau)$ and $P_{odd}(2 \tau) $ belong to $\Gamma_0(4)$.
Hence, in the same way as in Lemma \ref{lem:e1},  
we get our lemma.
\end{pf}

\section{moduli spaces on rational elliptic surfaces}

\subsection{Definition of invariants}

Let $ \pi:X \to {\Bbb P}^1$ be a rational elliptic surface
and $f$ a fiber of $\pi$.
$X$ is obtained by 9 points blow-ups of ${\Bbb P}^2$.
We assume that every fiber is irreducible. 
Let $H$ be the pull-back of the ample generator of 
$H^2({\Bbb P}^2,{\Bbb Z})$,
and $C_1,C_2,\dots,C_9$ the pull-backs of 
the classes of exceptional divisors.
Then  
$$
H^2(X,{\Bbb Z})={\Bbb Z} H \oplus {\Bbb Z} C_1
\oplus {\Bbb Z} C_2 \oplus \cdots \oplus {\Bbb Z} C_9.
$$
We set $g:=H-C_1$.
Then $g$ defines a morphism $\phi:X \to {\Bbb P}^1$,
which is factored to $X \to \Sigma_1 \to {\Bbb P}^1$.
For a sufficiently large integer $n$,
$nf+g$ becomes an ample divisor on $X$.
For a $c_1 \in H^2(X,{\Bbb Z})$ and a rational number $\Delta$,
 we can choose 
a divisor $h$ such that the equations for 
$\xi \in H^2(X,{\Bbb Z})$ with $-((2 \xi+c_1)^2) \leq 4 \Delta$
\begin{equation}
(2 \xi+c_1,f)=0,\quad
(2 \xi+c_1,g)=0,\quad
(2 \xi+c_1,h)=0
\end{equation}
have no other solution than $\xi=-c_1/2$.
For a polarization $L_{n_1,n_2}:=n_1 f+n_2 g+h$ with
$n_1 \gg n_2 \gg 1$,
we shall consider moduli spaces 
$M_{L_{n_1,n_2}}(c_1,\Delta)$.
It is known that $M_{L_{n_1,n_2}}(c_1,\Delta)$ is a 
smooth projective variety,
unless $2|c_1$ and $\Delta \in 2{\Bbb Z}$
[M].
We shall quot a useful lemma due to G\"{o}ttsche [G2].

\begin{lem}
Let $H$ be a divisor on a surface $X$ 
which belongs to the closure of the ample cone.
We assume that there is a non-negative number $k \in {\Bbb Q}$
such that $c_1(K_X)=-k c_1(H)$. 
For $c_1 \in H^2(X,{\Bbb Z})$ and $\Delta \in {\Bbb Q}$,
we shall choose a general ample divisor $L$ on $X$.
We assume that $M_{H+\epsilon L}(c_1,\Delta)$ is smooth.
Then $P(M_{H+\epsilon L}(c_1,\Delta),t)$, $0<\epsilon \ll 1$ 
does not depend on $L$ and $\epsilon$.
\end{lem} 
\begin{pf}
Since $\epsilon$ is sufficiently small,
every element $E$ of $M_{H+\epsilon L}(c_1,\Delta)$
is $\mu$-semi-stable with respect to $H$. 
Let $L'$ be another general ample divisor.
For rank 1 torsion free sheaves $I_1, I_2$ such that
$(c_1(I_1)-c_1(I_2),H)=0$,
$(c_1(I_1)-c_1(I_2),L')>0$ and
$(c_1(I_1)-c_1(I_2),L)<0$, we see that 
$\Ext^2(I_1,I_2) \cong \Hom(I_2,I_1 \otimes K_X)^{\vee}=0$.
Then we get that
\begin{align*}
& \sum_{\Delta}P(M_{H+\epsilon L}(l,\Delta),t)q^{\Delta}- 
\sum_{\Delta}P(M_{H+\epsilon' L'}(l,\Delta),t)q^{\Delta}\\
=&
\frac{1}{t(t-1)}
\left\{\sum
\begin{Sb}
\xi \in H^2(X,{\Bbb Z})+l/2\\
(\xi,H)=0\\
(\xi,L')> 0\\
(\xi,L)<0
\end{Sb}
 (t^2 q)^{-(\xi^2)}(t^{(\xi,K_X)}-
t^{-(\xi,K_X)})
\right\}=0.
\end{align*}
\end{pf}
Since we are only interested in Euler characteristics
of $M_{L_{n_1,n_2}}(c_1,\Delta)$,
we may denote $M_{L_{n_1,n_2}}(c_1,\Delta)$ by
$M(c_1,\Delta)$.
It is known that 
$$
\langle C_1,f \rangle^{\perp}
:=\{x \in H^2(X,{\Bbb Z})|(x,C_1)_*=(x,f)_*=0 \}
$$
is the $E_8$-lattice.
Hence we denote $\langle C_1,f \rangle^{\perp}$ by 
$E_8$.
We shall compute $\chi(M(c_1,\Delta)), c_1 \in E_8$.
The following proposition is due to Shioda [S].
\begin{prop}
There is a family of rational elliptic surfaces
$\pi_S: {\cal X} \to {\Bbb P}_S^1$ 
with a family of sections $\sigma:{\Bbb P}^1_S \to {\cal X}$
 parametrized by $S$ which satisfies that
\begin{enumerate}
\item
$W(E_8)$ acts on ${\cal X}$, $S$ and the action is compatible 
with the projection $p:{\cal X} \to S$.
\item
There is a $W(E_8)$-equivariant isomorphism
$$
\phi:E_8 \oplus {\Bbb Z}^{\oplus 2} \to R^2 p_* {\Bbb Z}
$$
such that $\phi({\Bbb Z}^{\oplus 2})=
{\Bbb Z}f \oplus {\Bbb Z} \sigma$ and
$\phi(E_8)=({\Bbb Z}f \oplus {\Bbb Z} \sigma)^{\perp}$,
where the action of $W(E_8)$ to ${\Bbb Z}^{\oplus 2}$ is trivial.
\end{enumerate}
\end{prop}
Since Maruyama [M] constructed a family of moduli of stable sheaves
on fibers ${\cal X}_s, s \in S$,
we obtain that,

\begin{prop}
We assume that $c_1$ belongs to 
$({\Bbb Z}f \oplus E_8) \setminus  2E_8$. 
For $w \in W(E_8)$,
$M_{L_{n_1,n_2}}(c_1,\Delta)$ is deformation equivalent
to $M_{w(L_{n_1,n_2})}(w(c_1),\Delta)$.
In particular, 
$P(M(c_1,\Delta),t)=P(M(w(c_1),\Delta),t)$, $w \in W(E_8)$.
\end{prop}
Since there are three orbits of the action of $W(E_8)$
to $E_8/2 E_8$, we define our invariants to compute
as follows (cf. [MNVW, (3.15)]).
\begin{defn}
For a class $c_1 \in H^2(X,{\Bbb Z})$,
we set
\begin{equation}
\widetilde{Z}_{c_1}(\tau):=
q^{-2 \chi(X)/24} \sum_{\Delta}
\chi(M(c_1,\Delta))q^{\Delta}.
\end{equation}

Let $v_0, v_{even}$ and $v_{odd}$ be elements of $E_8$
such that $v_0=0$, $(v_{even}^2)_*=4$ and 
$(v_{odd}^2)_*=2$.
We set
\begin{equation}
\widetilde{Z}_{\lambda}(\tau):=
\frac{1}{2}\left(\widetilde{Z}_{v_{\lambda}}(\tau)+
\widetilde{Z}_{f+v_{\lambda}}(\tau)\right),\;\lambda=0, even, odd.
\end{equation}
\end{defn}

\subsection{Poincar\'{e} polynomials of moduli spaces}

\begin{defn}
Let $l \in H^2(X,{\Bbb Z})$ be a cohomology class such that
$(l,g)=0$.\newline
$M_g(l,\Delta)$ is the set ( or the stack) of torsion free sheaves $E$
of rank 2 on $X$ such that
\begin{enumerate}
\item $E$ is semi-stable with respect to $g$, that is,
$E_{|D} \cong {\cal O}_D^{\oplus 2}$ for general member $D \in |g|$,
\item $(c_1(E),\Delta(E))=(l,\Delta)$.
\end{enumerate}
\end{defn}
We set
\begin{equation}
\begin{split}
F(t,q) &:=\frac{1}{(t^2-1)(t-1)}
\prod_{a \geq 1}Z_t(\Sigma_1,t^{2a-2}q^a)Z_t(\Sigma_1,t^{2a}q^a),\\
G(t,q) &:=\frac{1}{(t^2-1)(t-1)}
\left\{\prod_{a \geq 1}Z_t(\Sigma_1,t^{2a-2}q^a)Z_t(\Sigma_1,t^{2a}q^a)
-\prod_{a \geq 1}Z_t(\Sigma_1,t^{2a-1}q^a)^2 \right\},\\
B_0(t,u)& :=\sum_{n \in \Bbb Z}u^{n^2}t^n=\Theta_{A_1}(x,u), \\
B_1(t,u)& :=\sum_{n \in {\Bbb Z}+1/2}u^{n^2}t^n
=(\Theta_{A_1}|1/2)(x,u),
\end{split}
\end{equation}
where $e^x=t$ and ${A}_1$ is the ${A}_1$-lattice.

Motivated by [Y2, Thm. 0.1 and Lem. 1.4]
(also see [Y3, Cor. 3.5]), 
we define the Poincar\'{e} polynomial of $M_g(l,\Delta)$
( or $(t-1) \times \text{(the number of ${\Bbb F}_t$-rational points of
$M_g(l,\Delta)$)}$)
as follows:
\begin{equation}\label{eq:betti}
\sum_{\Delta}P(M_g(l,\Delta),t)q^{\Delta}:=
\frac{F(t,q) B_0(t,t^2q)^{8-n}B_1(t,t^2q)^n}
{\prod_{a \geq 1}(1-(t^2q)^a)^{16}},
\end{equation}
where $n=\#\{i|2 \leq i \leq 9,
(l,E_i) \equiv 1 \mod 2 \}$.

\begin{rem}
This definition is closely related to $SL$-equivariant
cohomologies of some spaces 
([Y2, sect. 4], [Y3, sect. 3]).
\end{rem}
 
In the same way as in the proof of
[Y2, Prop. 3.3] (or, [Y1, Rem. 4.6],
[G2]),
we see that
\begin{equation}\label{eq:wall}
\begin{split}
& \sum_{\Delta}P(M(l,\Delta),t)q^{\Delta}- 
\sum_{\Delta}P(M_g(l,\Delta),t)q^{\Delta}\\
=&
\frac{1}{t(t-1)}
\left\{\sum
\begin{Sb}
\xi \in H^2(X,{\Bbb Z})+l/2\\
(\xi,g)>0\\
(\xi,f) <0
\end{Sb}
 (t^2 q)^{-(\xi^2)}t^{(\xi,K_X)}-
\sum
\begin{Sb}
\xi \in H^2(X,{\Bbb Z})+l/2\\
(\xi,g) \leq 0\\
(\xi,f) > 0
\end{Sb}
 (t^2 q)^{-(\xi^2)}t^{(\xi,K_X)}\right\}
\prod_{a \geq 1}Z_t(X,t^{-1}(t^2q)^a)^2\\
& \quad -\frac{1}{2t(t-1)}\left\{
\sum
\begin{Sb}
\xi \in H^2(X,{\Bbb Z})+l/2\\
(\xi,g)=0\\
(\xi,f)=0
\end{Sb}
 (t^2 q)^{-(\xi^2)}t^{(\xi,K_X)}
\right\}\prod_{a \geq 1}Z_t(X,t^{-1}(t^2q)^a)^2.
\end{split}
\end{equation}
The last term comes from the contribution of
the set of $E$ which have the Harder-Narasimhan filtration
$$
0 \to I_1 \to E \to I_2 \to 0,
$$ 
such that $(c_1(I_1)-c_1(I_2),f)=(c_1(I_1)-c_1(I_2),g)=0$
and $(c_1(I_1)-c_1(I_2),h)>0$.

\begin{rem}\label{rem:def}
If $M(l,\Delta)$ is singular, that is,
$l=0$ and $\Delta$ is even,
then we shall define $\chi(M(l,\Delta))$
by this equation for a moment.
\end{rem}

The first term is closely related to Jacobi form of
G\"{o}ttsche and Zagier [G-Z1].
Hence we can expect that
$$
 \sum_{\Delta}P(M(l,\Delta),t)q^{\Delta}
$$ 
has a good expression. 
As in the proof of [G-Z1, Thm. 3.9],
we shall consider a sublattice
$\langle f,g \rangle \oplus  \langle f,g \rangle^{\perp} \subset
H^2(X,{\Bbb Z})$ and Jacobi forms associated with this lattice.
We set
\begin{equation}
\begin{cases}
e_i:=C_{9-i}-C_{10-i}, (1 \leq i \leq 7),\\
e_8:=H-C_1-C_2-C_3,\\
{\frak p}:=(e_1+e_3+e_5+e_8)/2,\\
{\frak q}:=(e_7+e_8)/2.
\end{cases}
\end{equation}
Then we see that
\begin{equation}
\begin{split}
\langle f,g \rangle^{\perp} &=
{\Bbb Z} e_1 \oplus {\Bbb Z} e_2 \oplus 
\cdots \oplus {\Bbb Z} e_8 \\
&= D_8,
\end{split}
\end{equation}
and 
\begin{equation}
H^2(X,{\Bbb Z})/
(\langle f,g \rangle \oplus \langle f,g \rangle^{\perp})
\cong {\Bbb Z}/2{\Bbb Z}\, ( f/2+{\frak p})
 \oplus {\Bbb Z}/2{\Bbb Z}\, (g/2+{\frak q}).
\end{equation}
We set
\begin{align*}
A_1(t,q) &:=\sum_{m,n>0}(t^2 q)^{4mn}\frac{t^{2m}-t^{-2m}}{t-1},\\
A_2(t,q) & :=\sum_{m,n>0}(t^2 q)^{2m(2n-1)}\frac{t^{2m}-t^{-2m}}{t-1},\\
A_3(t,q) & :=\sum_{m,n>0}(t^2 q)^{(2m-1)2n}\frac{t^{2m-1}-t^{-(2m-1)}}{t-1},\\
A_4(t,q) &:=\sum_{m,n>0}(t^2 q)^{(2m-1)(2n-1)}
\frac{t^{2m-1}-t^{-(2m-1)}}{t-1}.
\end{align*}
Then we get that
\begin{lem}\label{lem:j}
\begin{align*}
A_1(1,q) &=4\sum_{n>0}\sigma_1(n) q^{4n}=\frac{4(1-E_2(4 \tau))}{24},\\
A_2(1,q) & =4\sum_{m,n>0}m q^{2m(2n-1)}=
\frac{4(E_2(4\tau)-E_2(2 \tau))}{24},\\
A_3(1,q) & =2\sum_{n>0}\sigma_1^{odd}(n)q^{2n}=
\frac{-E_2(2\tau)+2E_2(4 \tau)-1}{12},\\
A_4(1,q) & =\sum_{m,n>0} 2(2m-1) q^{(2m-1)(2n-1)}=2 F(\tau).
\end{align*}
\end{lem}
The following is an easy consequence of the description of $D_8$
in ${\Bbb Z}^8$ [C-S, Chap. 4]. 
\begin{lem}\label{lem:d8}
\begin{align*}
\Theta_{D_8}(x,u) &=
\frac{\Theta_{\Bbb Z}(x,u)^8+\Theta_{\Bbb Z}(x+1/2,u)^8}{2},\\
(\Theta_{D_8}|{\frak q})(x,u) &=
\frac{\Theta_{\Bbb Z}(x,u)^8-\Theta_{\Bbb Z}(x+1/2,u)^8}{2},\\
(\Theta_{D_8}|{\frak p})(x,u) &=
(\Theta_{D_8}|({\frak p}+{\frak q}))(x,u)\\
&=\frac{(\Theta_{\Bbb Z}|(1/2))(x,u)^8}{2},\\
(\Theta_{D_8}|(e_1/2))(x,u) &=(\Theta_{D_8}|(e_1/2+{\frak q}))(x,u)\\ 
&=\frac{\Theta_{\Bbb Z}(x,u)^6 (\Theta_{\Bbb Z}|(1/2))(x,u)^2}{2},\\
(\Theta_{D_8}|(e_1/2+{\frak p}))(x,u) &
=(\Theta_{D_8}|(e_1/2+{\frak p}+{\frak q}))(x,u)\\ 
&=\frac{\Theta_{\Bbb Z}(x,u)^2 (\Theta_{\Bbb Z}|(1/2))(x,u)^6}{2}.
\end{align*}
\end{lem}

\section{Computation of $\widetilde{Z}_{v_{\lambda}}(\tau)$}

\subsection{Computation of $\widetilde{Z}_{v_{even}}(\tau)$}

We first treat the case of $Z_{v_{even}}$.
Hence we assume that $c_1=e_7+e_8$.
By using \eqref{eq:betti}, \eqref{eq:wall} and the proof of
[G-Z1, Thm. 3.9],
we obtain the following proposition.

\begin{prop}
\begin{multline}
\sum_{\Delta}P(M(e_7+e_8,\Delta),t)q^{\Delta}=
\frac{F(t,q) B_0(t,u)^8}{\prod_{a \geq 1}(1-u^a)^{16}} \\
+\prod_{a \geq 1}Z_t(X,t^{-1}u^a)^2
\left\{-\frac{1}{t(t^2-1)(t-1)} (\Theta_{D_8}|{\frak q})(0,u^2)
-\frac{1}{(t^2-1)(t-1)} \Theta_{D_8}(0,u^2) \right.\\
+\frac{A_1(t,q) (\Theta_{D_8}|{\frak q})(0,u^2)}{t}\\
\left.
+\frac{A_2(t,q) (\Theta_{D_8}|({\frak p}+{\frak q}))(0,u^2)}{t}
+\frac{A_3(t,q)\Theta_{D_8}(0,u^2)}{t} \right.\\
\left. +\frac{A_4(t,q)(\Theta_{D_8}|{\frak p})(0,u^2)}{t}   
-
\frac{(\Theta_{D_8}|{\frak q})(0,u^2)}
{2t(t-1)}\right\},
\end{multline}
where $u=u(t):=t^2 q$.
\end{prop}
By Lemma \ref{lem:d8},
$(\Theta_{D_8}+\Theta_{D_8}|{\frak q})(0,u^2)=B_0(1,u)^8$.
Hence we get that
\begin{align*}
\sum_{\Delta}P(M(e_7+e_8,\Delta),t)q^{\Delta}
=&\frac{G(t,q) B_0(t,u)^8}{\prod_{a \geq 1}(1-u^a)^{16}}
+\prod_{a \geq 1}Z_t(X,t^{-1}u^a)^2
\left\{\frac{B_0(t,u)^8-B_0(1,u)^8
}{(t+1)(t-1)^2}\right. \\
&+\left(\frac{1}{(t^2-1)(t-1)}-\frac{1}{t(t^2-1)(t-1)}\right)
 (\Theta_{D_8}|{\frak q})(0,u^2)\\
&+\frac{A_1(t,q) (\Theta_{D_8}|{\frak q})(0,u^2)}{t}
+\frac{A_2(t,q) (\Theta_{D_8}|({\frak p}+{\frak q}))(0,u^2)}{t}\\
&+\frac{A_3(t,q)\Theta_{D_8}(0,u^2)}{t}\\
& \left. +\frac{A_4(t,q)(\Theta_{D_8}|{\frak p})(0,u^2)}{t}   
-
\frac{(\Theta_{D_8}|{\frak q})(0,u^2)}
{2t(t-1)}\right\}\\
= &  \frac{G(t,q) B_0(t,u)^8}{\prod_{a \geq 1}(1-u^a)^{16}}
+\prod_{a \geq 1}Z_t(X,t^{-1}u^a)^2
\left\{\frac{B_0(t,u)^8-B_0(1,u)^8
}{(t+1)(t-1)^2}\right.  \\
&
-\frac{1}{2t(t+1)}
 (\Theta_{D_8}|{\frak q})(0,u^2)
+\frac{A_1(t,q) (\Theta_{D_8}|{\frak q})(0,u^2)}{t}\\
&\left.
+\frac{A_2(t,q) (\Theta_{D_8}|({\frak p}+{\frak q}))(0,u^2)}{t}
+\frac{A_3(t,q)\Theta_{D_8}(0,u^2)}{t} \right.\\
&\left. +\frac{A_4(t,q)(\Theta_{D_8}|{\frak p})(0,u^2)}{t}
\right\},
\end{align*}
where $u=t^2 q$.
In order to evaluate at $t=1$,
we need the following two lemmas.

\begin{lem}\label{lem:ruled}
\begin{equation}
q^{-1/3}G(1,q)=\frac{2\sum_{n \geq 1} \sigma_1(n)q^n}{\eta(q)^8}
=\frac{1-E_2(\tau)}{12 \eta(\tau)^8}.
\end{equation}

\end{lem}

\begin{lem}\label{lem:blow-up}
\begin{equation}
\begin{split}
\lim_{t \to 1}\frac{B_0(t,t^2q)-B_0(1,t^2q)}{(t+1)(t-1)^2}
&=-\frac{1}{2}\left(\sum_{m,n>0}(-1)^m m q^{m(2n-1)} \right)B_0(1,q)\\
\lim_{t \to 1}\frac{B_1(t,t^2q)-B_1(1,t^2q)}{(t+1)(t-1)^2}
&=-\frac{1}{2}\left(\sum_{m,n>0}(-1)^m m q^{2mn}-\frac{1}{8} \right)B_1(1,q).
\end{split}
\end{equation}
\end{lem}

\begin{pf}
By the product formula,
\begin{equation}
B_0(t,u)=
\prod_{m \geq 1}(1-u^{2m})(1+u^{2m-1}t)(1+u^{2m-1}t^{-1}).
\end{equation}
Then we see that
\begin{equation}
\frac{d \log B_0}{d t}=
\sum_{m>0} \frac{u^{2m-1}(t^2-1)}
{t^2(1+u^{2m-1}t)(1+u^{2m-1}t^{-1})}.
\end{equation}
Hence we see that
\begin{align*}
\lim_{t \to 1}\frac{B_0(t,u)-B_0(1,u)}{(t-1)^2}
&=\frac{1}{2}\frac{d^2 B_0}{d t^2}(1,u)\\
&=\sum_{m \geq 1} \frac{u^{2m-1}}{(1+u^{2m-1})^2}\\
&=-\sum_{m,n>0}(-1)^n n u^{(2m-1)n}.
\end{align*}
Thus we obtain the first relation.
For the second relation,
we use the identity
\begin{equation}
B_1(t,u)=
u^{1/4}t^{1/2}\prod_{m \geq 1}(1-u^{2m})(1+u^{2m}t)(1+u^{2m-2}t^{-1}).
\end{equation}
\end{pf}
Combining Lemma \ref{lem:wt2}, \ref{lem:j}, \ref{lem:d8}, \ref{lem:ruled},
\ref{lem:blow-up}
all together, we see that
\begin{align*}
\widetilde{Z}_{v_{even}}(\tau)
&=\lim_{t \to 1}q^{-1} \sum_{\Delta}
 P(M(e_7+e_8,\Delta),t)q^{\Delta}\\
&= \frac{1}{24\eta(\tau)^{24}}
\left(\theta_3(2 \tau)^8(-6 E_2(\tau)+19 E_2(2\tau)-16 E_2(4\tau)) \right.\\
& \quad \left. -\theta_4(2 \tau)^8(-E_2(2 \tau)+4 E_2(4 \tau))
         -\theta_2(2 \tau)^8(E_2(\tau)-E_2(2 \tau))\right).
\end{align*}
The following is an easy consequence of Lemma \ref{lem:e1}.
\begin{lem}
\begin{equation}
\begin{split}
& 15(\theta_3(2 \tau)^8-\theta_4(2 \tau)^8)e_1(\tau)
+48(4\theta_3(2 \tau)^8-\theta_4(2 \tau)^8)F(\tau)
-3 \theta_2(2\tau)^8 e_1(\tau)\\
= &
\frac{1}{8}  \theta_2(\tau)^8(\theta_3(\tau)^4+\theta_4(\tau)^4).
\end{split}
\end{equation}

\end{lem}
Combining all together, we obtain that
\begin{equation}\label{eq:main1}
\begin{split}
\widetilde{Z}_{v_{even}}(\tau)=&
 \frac{1}{24\eta(\tau)^{24}}
\left(\theta_3(2 \tau)^8(-6 E_2(\tau)+19 E_2(2\tau)-16 E_2(4\tau)) \right.\\
& \quad \left. +\theta_4(2 \tau)^8(-E_2(2 \tau)+4 E_2(4 \tau))
         -\theta_2(2 \tau)^8(E_2(\tau)-E_2(2 \tau))\right)\\
=& \frac{1}{24\eta(\tau)^{24}}
\left(-\frac{\theta_3(2 \tau)^8-\theta_4(2 \tau)^8+\theta_2(2 \tau)^8}{2}
E_2(\tau)+\theta_3(2\tau)^8(15 e_1(\tau)+4 \cdot 48 F(\tau)) \right.\\
& \quad -\left. \theta_4(2\tau)^8(15 e_1(\tau)+48F(\tau))
-3\theta_2(2 \tau)^8 e_1(\tau)\right )\\
=& -\frac{1}{24 \eta(\tau)^{24}}
\left(\frac{E_2(\tau)P_{even}(\tau)}{135}-
\theta_2(\tau)^8 \frac{(\theta_3(\tau)^4+\theta_4(\tau)^4)}{8}\right).
 \end{split}
\end{equation}

\subsection{Computation of $\widetilde{Z}_{v_{0}}(\tau)$}

By using \eqref{eq:betti}, \eqref{eq:wall} and the proof of
[G-Z1, Thm. 3.9] again, we get that
\begin{prop}
\begin{multline}
\sum_{\Delta}P(M(0,\Delta),t)q^{\Delta}=
\frac{F(t,q) B_0(t,u)^8}{\prod_{a \geq 1}(1-u^a)^{16}} \\
+\prod_{a \geq 1}Z_t(X,t^{-1}u^a)^2
\left\{-\frac{1}{t(t^2-1)(t-1)} \Theta_{D_8}(0,u^2)
-\frac{1}{(t^2-1)(t-1)} (\Theta_{D_8}|{\frak q})(0,u^2)\right.\\
+\frac{A_1(t,q) \Theta_{D_8}(0,u^2)}{t} 
+\frac{A_2(t,q) (\Theta_{D_8}|{\frak p})(0,u^2)}{t}
+\frac{A_3(t,q)(\Theta_{D_8}|{\frak q})(0,u^2)}{t} \\
\left. +
\frac{A_4(t,q)(\Theta_{D_8}|({\frak p}+{\frak q}))(0,u^2)}{t}   
-
\frac{\Theta_{D_8}(0,u^2)}
{2t(t-1)}\right\},
\end{multline}
where $u=u(t):=t^2 q$.
\end{prop}
In the same way as in the even case, we see that
\begin{align*}
\widetilde{Z}_{v_{0}}(\tau)
&=\lim_{t \to 1}q^{-1} \sum_{\Delta}
 P(M(0,\Delta),t)q^{\Delta}\\
&= \frac{1}{24\eta(\tau)^{24}}
\left(\theta_3(2 \tau)^8(-6 E_2(\tau)+19 E_2(2\tau)-16 E_2(4\tau)) \right.\\
& \quad \left. +\theta_4(2 \tau)^8(-E_2(2 \tau)+4 E_2(4 \tau))
         -\theta_2(2 \tau)^8(E_2(\tau)-E_2(2 \tau))\right).
\end{align*}
By using Lemma \ref{lem:e1} and Corollary \ref{cor:t1},
we see that
\begin{equation}
2 \theta_4(2 \tau)^8
(15 e_1(\tau)+3 \cdot 16 F(\tau))=
-\left(\theta_3(\tau)^4 \theta_4(\tau)^4(\theta_3(\tau)^4+
\theta_4(\tau)^4)+3 \theta_4(2\tau)^{12} \right).
\end{equation}
Comparing $\widetilde{Z}_{v_{0}}(\tau)$ with
$\widetilde{Z}_{v_{even}}(\tau)$, we obtain that

\begin{equation}\label{eq:main2}
\begin{split}
\widetilde{Z}_{v_{0}}(\tau)=&
 \frac{1}{24\eta(\tau)^{24}}
\left(\theta_3(2 \tau)^8(-6 E_2(\tau)+19 E_2(2\tau)-16 E_2(4\tau)) \right.\\
& \quad \left. +\theta_4(2 \tau)^8(-E_2(2 \tau)+4 E_2(4 \tau))
         -\theta_2(2 \tau)^8(E_2(\tau)-E_2(2 \tau))\right)\\
=& \frac{1}{24\eta(\tau)^{24}}
\left(-\frac{\theta_3(2 \tau)^8+\theta_4(2 \tau)^8+\theta_2(2 \tau)^8}{2}
E_2(\tau)+\theta_3(2\tau)^8(15 e_1(\tau)+4 \cdot 48 F(\tau)) \right.\\
& \quad +\left. \theta_4(2\tau)^8(15 e_1(\tau)+48F(\tau))
-3\theta_2(2 \tau)^8 e_1(\tau)\right )\\
=& \frac{-1}{24 \eta(\tau)^{24}}
\left(E_2(\tau)P_{0}(\tau)+
\left(\theta_3(\tau)^4\theta_4(\tau)^4-
\frac{\theta_2(\tau)^8}{8}\right)
(\theta_3(\tau)^4+\theta_4(\tau)^4)+
3\theta_4(2 \tau)^{12}\right)\\
=& \frac{-1}{24 \eta(\tau)^{24}}
\left(E_2(\tau)P_{0}(\tau)+
\left(\theta_3(\tau)^4\theta_4(\tau)^4-
\frac{\theta_2(\tau)^8}{8}\right)
(\theta_3(\tau)^4+\theta_4(\tau)^4)\right)-
\frac{1}{8 \eta(2 \tau)^{12}}.
\end{split}
\end{equation}

\begin{rem}\label{rem:euler}
Let $\widetilde{Z}_{v_{0}}^{int}(\tau)$ be the generating 
function of Euler characteristics associated to the intersection 
cohomologies $IH^*(M(0,n),{\Bbb Q})$.
Then in the same way as in [Y2, Rem. 4.6],
we get that
\begin{equation}
\widetilde{Z}_{v_{0}}^{int}(\tau)=
\widetilde{Z}_{v_{0}}(\tau)+\frac{1}{4 \eta(2 \tau)^{12}},
\end{equation}
where we used the definition of $\widetilde{Z}_{v_{0}}(\tau)$
in Remark \ref{rem:def}.
\end{rem}

\subsection{Computation of $\widetilde{Z}_{v_{odd}}(\tau)$}

Since $(e_1^2)_*=2$, we may assume that $v_{odd}=e_1=C_8-C_9$.
Then similar computations show that
\begin{prop}
\begin{multline}
\sum_{\Delta}P(M(e_1,\Delta),t)q^{\Delta}=
\frac{F(t,q) B_0(t,u)^6 B_1(t,u)^2}
{\prod_{a \geq 1}(1-u^a)^{16}} \\
+\prod_{a \geq 1}Z_t(X,t^{-1}u^a)^2
\left\{-\frac{1}{t(t^2-1)(t-1)} (\Theta_{D_8}|(e_1/2))(0,u^2)
-\frac{1}{(t^2-1)(t-1)} (\Theta_{D_8}|(e_1/2+{\frak q}))(0,u^2)\right.\\
+\frac{A_1(t,q) (\Theta_{D_8}|(e_1/2))(0,u^2)}{t} 
+\frac{A_2(t,q) (\Theta_{D_8}|(e_1/2+{\frak p}))(0,u^2)}{t}
+\frac{A_3(t,q)(\Theta_{D_8}|(e_1/2+{\frak q}))(0,u^2)}{t}\\
\left. 
+\frac{A_4(t,q)(\Theta_{D_8}|(e_1/2+{\frak p}+{\frak q}))(0,u^2)}{t}   
-
\frac{(\Theta_{D_8}|(e_1/2))(0,u^2)}
{2t(t-1)}\right\},
\end{multline}
where $u=u(t):=t^2 q$.
\end{prop}
Hence we get that
\begin{equation}\label{eq:main3}
\begin{split}
\widetilde{Z}_{v_{odd}}(\tau)
&=\lim_{t \to 1}q^{-1} \sum_{\Delta}
 P(M(e_1,\Delta),t)q^{\Delta}\\
&= \frac{1}{24\eta(\tau)^{24}}
\left(\theta_2(2 \tau)^2 \theta_3(2 \tau)^6
(-5 E_2(\tau)+13 E_2(2\tau)-8 E_2(4\tau)) \right.\\
& \quad \left. +\theta_2(2 \tau)^6 \theta_3(2 \tau)^2
(-E_2(\tau)+E_2(2 \tau))\right)\\
&=\frac{1}{24\eta(\tau)^{24}}
\left(-\frac{\theta_2(2 \tau)^2 \theta_3(2 \tau)^6+
\theta_2(2 \tau)^6 \theta_3(2 \tau)^2}{2}E_2(\tau)\right.\\
&\left. +\theta_2(2 \tau)^2 \theta_3(2 \tau)^6(9\cdot16 F(\tau)-3 e_1(\tau))
+\theta_2(2 \tau)^6 \theta_3(2 \tau)^2(-3 e_1(\tau))\right)\\
&=- \frac{1}{24\eta(\tau)^{24}}
\left(E_2(\tau) \frac{P_{odd}(\tau)}{120}-\frac{1}{8}\theta_2(\tau)^4 E_4(\tau)
\right).
\end{split}
\end{equation}

\subsection{Main results}
In the same way as in section 3.3, we can compute
$\widetilde{Z}_{f+v_{\lambda}}(\tau)$.
In particular, we obtain the following relations:
\begin{align*}
\widetilde{Z}_{f+v_0}(\tau) &=
\widetilde{Z}_{v_0}(\tau)+\frac{1}{4 \eta(2 \tau)^{12}},\\
\widetilde{Z}_{f+v_{even}}(\tau) &=
\widetilde{Z}_{v_{even}}(\tau),\\
\widetilde{Z}_{f+v_{odd}}(\tau) &=
\widetilde{Z}_{v_{odd}}(\tau).
\end{align*}

By using \eqref{eq:main1}, \eqref{eq:main2} and \eqref{eq:main3},
we obtain the following theorem. 
\begin{thm}\label{thm:1}
\begin{align*}
\widetilde{Z}_0(\tau)&=\frac{-1}{24 \eta(\tau)^{24}}
\left(E_2(\tau)P_{0}(\tau)+
\left(\theta_3(\tau)^4\theta_4(\tau)^4-
\frac{\theta_2(\tau)^8}{8}\right)
(\theta_3(\tau)^4+\theta_4(\tau)^4)\right),\\
\widetilde{Z}_{even}(\tau)&=\frac{-1}{24 \eta(\tau)^{24}}
\left(E_2(\tau)\frac{P_{even}(\tau)}{135}-
\theta_2(\tau)^8 \frac{(\theta_3(\tau)^4+\theta_4(\tau)^4)}{8}\right),\\
\widetilde{Z}_{odd}(\tau)&= \frac{-1}{24\eta(\tau)^{24}}
\left(E_2(\tau) \frac{P_{odd}(\tau)}{120}
-\frac{1}{8}\theta_2(\tau)^4 E_4(\tau)
\right).
\end{align*}
In particular,
$\widetilde{Z}_{\lambda}(\tau)=
Z_{\lambda}(\tau)$, up to holomorphic anomaly coming from 
$E_2(\tau)$ and unknown terms coming from singular spaces
$M(0,2n)$.
\end{thm}

By Remark \ref{rem:euler} and our theorem,
we conjecture the following.
\begin{guess}
For rational surfaces,
\begin{equation}\label{eq:guess}
Z_{v_{0}}(\tau)=\widetilde{Z}_{v_{0}}^{int}(\tau)
-\frac{1}{4 \eta(2 \tau)^{\chi(X)}},
\end{equation}
up to holomorphic anomaly.
\end{guess}

As a final remark, we shall check the S-duality conjecture in [V-W].
We set 
$w_i:=C_1+if, i=0,1$.
In [Y4], we showed that
\begin{align*}
Z_{w_0}&=\frac{1}{2}
\left(\frac{1}{\eta(\tau/2)^{12}}+\frac{\sqrt{-1}}
{\eta(\tau/2+1/2)^{12}}\right),\\
Z_{w_1}&=\frac{1}{2}
\left(\frac{1}{\eta(\tau/2)^{12}}-\frac{\sqrt{-1}}
{\eta(\tau/2+1/2)^{12}}\right).
\end{align*}
Hence we obtain that
\begin{align*}
Z_{SU(2)}(\tau) &=\frac{1}{2}
\left({Z}_0(\tau)-\frac{1}{8 \eta(2 \tau)^{12}} \right),\\
Z_{SO(3)}(\tau)&=2({Z}_0(\tau)+135{Z}_{even}(\tau)
+120{Z}_{odd}(\tau))+\frac{256}{\eta(\tau/2)^{12}}.
\end{align*}
Then we see that
\begin{equation}
Z_{SU(2)}(-1/\tau)=-2^{-6}\left(\frac{\tau}{\sqrt{-1}}\right)^{-6}
Z_{SO(3)}(\tau).
\end{equation}
Thus the S-duality conjecture in [V-W] hold for this case.

{\it 
Acknowledgement.}
This work is motivated by [MNVW] (especially by 
the form of (6.15)) and valuable discussions
with C. Vafa.
In particular, modular forms $\widehat{E}_2$ and $\eta^{24}$ 
in (6.15) made me to use [Y2].
I would like to thank J.A. Minahan, D. Nemeschansky, 
N.P. Warner, and especially 
C. Vafa very much.
I would also like to thank L. G\"{o}ttsche, H. Kanno, M-H. Saito, and Y. Yamada
for valuable discussions.

\end{document}